\documentclass[12pt]{article}
  \usepackage[cp1251]{inputenc}
  \usepackage[english]{babel}
  \usepackage{amsfonts,amssymb,amsmath,amstext,amsbsy,amsopn,amscd,graphicx}
  \usepackage{graphics}
  \usepackage{amsthm}

 \newtheorem{theorem}{Theorem}[section]
 \newtheorem{lemma}{Lemma}[section]
 \newtheorem{crl}{Corollary}[section]

 \def\thtext#1{
 \catcode`@=11
 \gdef\@thmcountersep{. #1}
 \catcode`@=12}

 \newcounter{Rk}[section]
 \renewcommand{\thtext}{\thesection.\arabic{Rk}}
 \newenvironment{remark}{\trivlist\item[\hskip\labelsep{\bf Remark}]
 \refstepcounter{Rk}{\bf\thesection.\arabic{Rk}.}}%
 {\endtrivlist}

 \newcounter{Df}[section]
 \renewcommand{\thtext}{\thesection.\arabic{Df}}
 \newenvironment{definition}{\trivlist\item[\hskip\labelsep{\bf Definition}]
 \par\refstepcounter{Df}{\bf\thesection.\arabic{Df}.}}%
 {\endtrivlist}

 \newenvironment{example}{\trivlist \item[\hskip\labelsep{\bf Example.}]}%
 {\endtrivlist}

\def\Z{\mathbb{Z}}

\def\0{{\mathbf 0}}
\def\b1{{\mathbf 1}}

\def\bb{{\mathbf b}}

\def\sign{\operatorname{sign}}
\def\corank{\operatorname{corank}}

 \title {An equivalence between the set of graph-knots and the set of homotopy classes of looped graphs}

\author{D.\,P.~Ilyutko\footnote{Partially supported by grants of RF President
NSh -- 660.2008.1, RFBR 07--01--00648, RNP 2.1.1.3704, the Federal
Agency for Education NK-421P/108.}}
\date {}

\begin {document}

 \maketitle
 \abstract{In the present paper we construct a one-to-one
 correspondence between the set of graph-knots and the set of
 homotopy classes of looped graphs. Moreover, the graph-knot and the
 homotopy class constructed from a given knot are related with this
 correspondence. This correspondence is given by a simple formula.}


 \section {Introduction}


The discovery of virtual knot theory by Kauffman~\cite{KaV} in mid
1990s was an important step in generalising combinatorial and
topological knot theoretical techniques into a larger domain (knots
in thickened surfaces), which is an important step towards
generalisation of these techniques for knots in arbitrary manifolds.

It turned out that some invariants (Kauffman bracket polynomial) can
be generalised for virtual knots immediately~\cite{KaV}, and some
other theories (Khovanov homology theory) need a complete revision
of the original construction for a generalisation for the case of
virtual knots~\cite{Izv}. Virtual knots also sharpened several
problems and elicited some phenomena which do not appear in the
classical knot theory~\cite{FKM}, e.g., the existence of a virtual
knot with non-trivial Jones polynomial and trivial fundamental group
emphasises the difficulty of extracting the Jones polynomial
information out of the knot group.

In the present paper, we consider two new theories: {\em
graph-links} introduced in~\cite {IM1,IM2} and {\em homotopy classes
of looped interlacement graphs} introduced by L.~Traldi and
L.~Zulli~\cite {TZ}. The two theories are closely connected to both
classical and virtual knot theories and, in some sense, are
generalisations of virtual knots. Likewise virtual knots appear out
of non-realisable Gauss code and thus generalise classical knots
(which have realisable Gauss codes), graphs-links and looped graphs
come out of {\em intersection graphs}: we may consider graphs which
realise chord diagrams, and, in turn, virtual knots, and pass to
arbitrary graphs which correspond to some mysterious objects
generalising knots and virtual knots.

 \begin{figure} \centering\includegraphics[width=200pt]{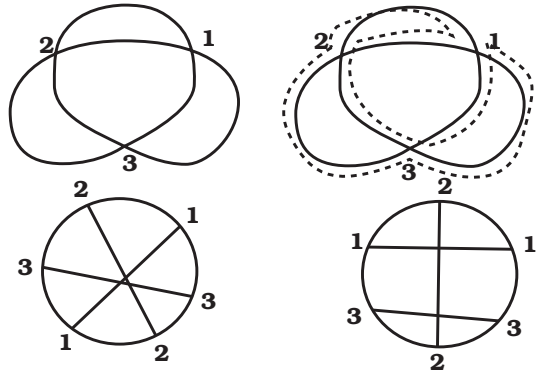}
\caption{A Gauss Circuit; A Rotating Circuit} \label{Gaussrot}
 \end{figure}

There is a way of coding all virtual knots by {\em Gauss diagrams}
and there is another way of coding all virtual links by {\em
rotating circuits} (see, e.g.,~\cite{Arxiv}), see Fig.~\ref
{Gaussrot}. Looped graphs are a generalisation of intersection
graphs of chord diagrams constructed by using the first way of
coding and graph-links come out from considering a rotating circuit.
In low-dimensional topology both approaches, the Gauss circuit
approach and the rotating circuit approach, are very widely used.
The Gauss circuit approach is applied in knot theory, namely in the
construction of finite-type invariants, Vassiliev invariants~\cite
{BN,GPV,CDL}, and in the planarity problem of immersed curves,~\cite
{CE1,CE2,RR}. However, for detecting planarity of a framed 4-graph
it is more convenient to use the rotating circuit approach,
see~\cite {Ma1,Ma2}. The criterion of the planarity of an immersed
curve, which is the framed 4-graph, is formulated very easy: an
immersed curve is planar if and only if the chord diagram obtained
from a rotating circuit is a {\em d-diagram}, i.e.\ the set of all
the chords can be split into two sets and the chords from one set do
not intersect each other, see~\cite {Ma3}. If we want to generalise
the planarity problem to the problem of finding the minimum genus of
a closed surface which a given curve can be immersed in, the
rotating circuit approach is also more useful. There are criteria
giving us the answer to the question what is the minimum genus for a
given curve, see~\cite {Ma1}.

 \begin{figure}
\centering\includegraphics[width=300pt]{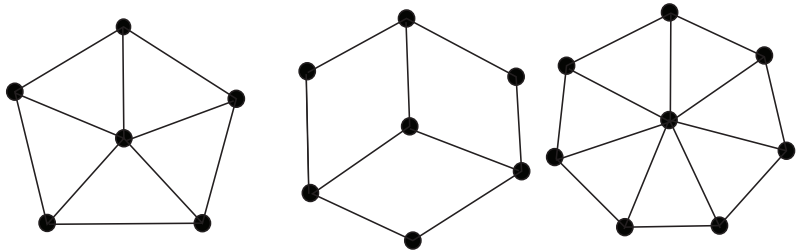}
\caption{Non-Realisable Bouchet Graphs} \label{Bouchet}
 \end{figure}

In spite of the fact that graph-links and the homotopy classes of
looped graphs are abstract objects the Kauffman bracket polynomial
and the Jones polynomial were constructed for them, see~\cite
{IM1,IM2,TZ}. But the first question which has arisen after the
constructions of these two theories was whether or not in each
graph-link (each homotopy class of looped graphs) there exists a
`realisable' representative, i.e.\ a graph which is the intersection
graph of a chord diagram. Some graphs can not be represented by
chord diagrams at all~\cite {Bou}, see Fig.~\ref{Bouchet}. The
answer to this question is negative and this problem was first
solved by V.\,O.~Manturov in~\cite {Ma4,Ma5} for homotopy classes of
looped graphs. The second question is whether there is an
equivalence between the two theories. In the present paper we give
an explicit formula which gives us an equivalence between the set of
graph-knots (graph-links with one component) and the set of homotopy
classes of looped graphs. Moreover, under this equivalence
`realisable' objects correspond to `realisable' ones. Therefore, by
using this equivalence we can give the negative answer to the
question of existing a realisable representative for graph-knots.

 \section*{Acknowledgments}
The author is grateful to V.\,O.~Manturov, A.\,T.~Fomenko,
D.\,M.~Afa\-na\-siev, I.\,M.~Nikonov for their interest to this work
and to L.~Traldi for his discussion of~\cite {TZ} and useful
comments.

 \section{Basic Definitions and Constructions}

\subsection{Chord diagrams and Framed 4-Graphs}

Throughout the paper all graphs are finite. Let $G$ be a graph with
the set of vertices $V(G)$ and the set of edges $E(G)$. We think of
an edge as an equivalence class of the two half-edges. We say that a
vertex $v\in V(G)$ has {\em degree} $k$ if $v$ is incident to $k$
half-edges. A graph whose vertices have the same degree $k$ is
called {\em $k$-valent} or a {\em k-graph}. The free loop, i.e.\ the
graph without vertices, is also considered as $k$-graph for any $k$.

 \begin {definition}\label {def:fr4}
A $4$-graph is {\em framed} if for every vertex the four emanating
half-edges are split into two pairs of (formally) opposite edges.
The edges from one pair are called {\em opposite to each other}.
 \end {definition}

\begin {figure} \centering\includegraphics[width=300pt]{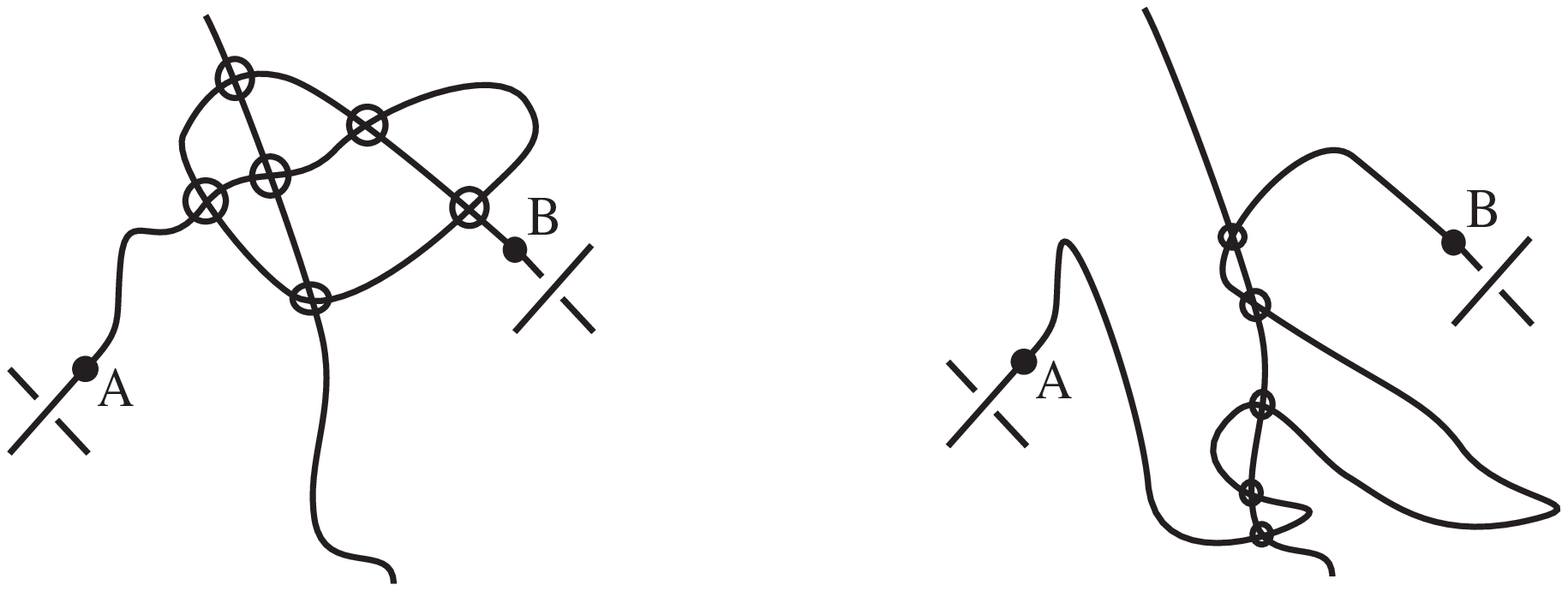}
  \caption{The detour move} \label{detour}
 \end {figure}

A {\em virtual diagram} is a framed $4$-graph embedded into
${\mathbb R}^2$ where each crossing is either endowed with a
classical crossing structure (with a choice for underpass and
overpass specified) or just said to be virtual and marked by a
circle. A {\em virtual link} is an equivalence class of virtual
diagrams modulo generalised Reidemeister moves. The latter consist
of usual Reidemeister moves referring to classical crossings and the
{\em detour move} that replaces one arc containing only virtual
intersections and self-intersections by another arc of such sort in
any other place of the plane, see Fig.~\ref{detour}. A {\em
projection} of a virtual diagram is a framed $4$-graph obtained from
the diagram by considering classical crossings as vertices and
virtual crossings are just intersection points of images of
different edges. A virtual diagram is {\em connected} if its
projection is connected. Without loss of generality, {\it all
virtual diagrams are assumed to be connected and contain at least
one classical crossing}~\cite {IM1,IM2}.

  \begin{definition}
A {\em chord diagram} is a cubic graph consisting of a selected
cycle (the {\em circle}) and several non-oriented edges ({\em
chords}) connecting points on the circle in such a way that every
point on the circle is incident to at most one chord. A chord
diagram is {\em labeled} if every chord is endowed with a label
$(a,\alpha)$, where $a\in\{0,1\}$ is the framing of the chord, and
$\alpha\in\{\pm\}$ is the sign of the chord. If no labels are
indicated, we assume the chord diagram has all chords with label
$(0,+)$. Two chords of a chord diagram are called {\em linked} if
the ends of one chord lie in different connected components of the
circle with the end-points of the second chord removed.
  \end{definition}

\begin{figure} \centering\includegraphics[width=200pt]{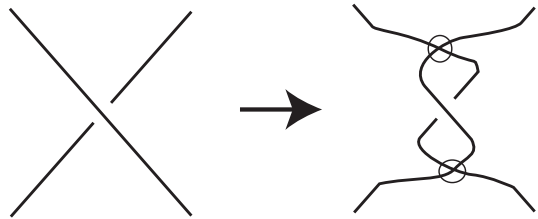}
  \caption{Virtualisation} \label{virtualisation}
 \end{figure}

 \begin {definition}
By a {\em virtualisation} of a classical crossing of a virtual
diagram we mean a local transformation shown in
Fig.~\ref{virtualisation}.
 \end {definition}

Having a labeled chord diagram $D$, one can construct a virtual link
diagram $K(D)$ (up to virtualisation) as follows. Let us immerse
this diagram in ${\mathbb R}^{2}$ by taking an embedding of the
circle and placing some chords inside the circle and the other ones
outside the circle. After that we remove neighbourhoods of each of
the chord ends and replace them by a pair of lines (connecting four
points on the circle which are obtained after removing two
neighbourhoods) with a classical crossing if the chord is framed by
$0$ and a couple of lines with a classical crossing and a virtual
crossing if the chord is framed by $1$ in the following way. The
choice for underpass and overpass is specified as follows. A
crossing can be smoothed in two ways: $A$ and $B$ as in the Kauffman
bracket polynomial; we require that the initial piece of the circle
corresponds to the $A$-smoothing if the chord is positive and to the
$B$-smoothing if it is negative:
 $A\colon\raisebox{-0.25\height}{\includegraphics[width=0.5cm]{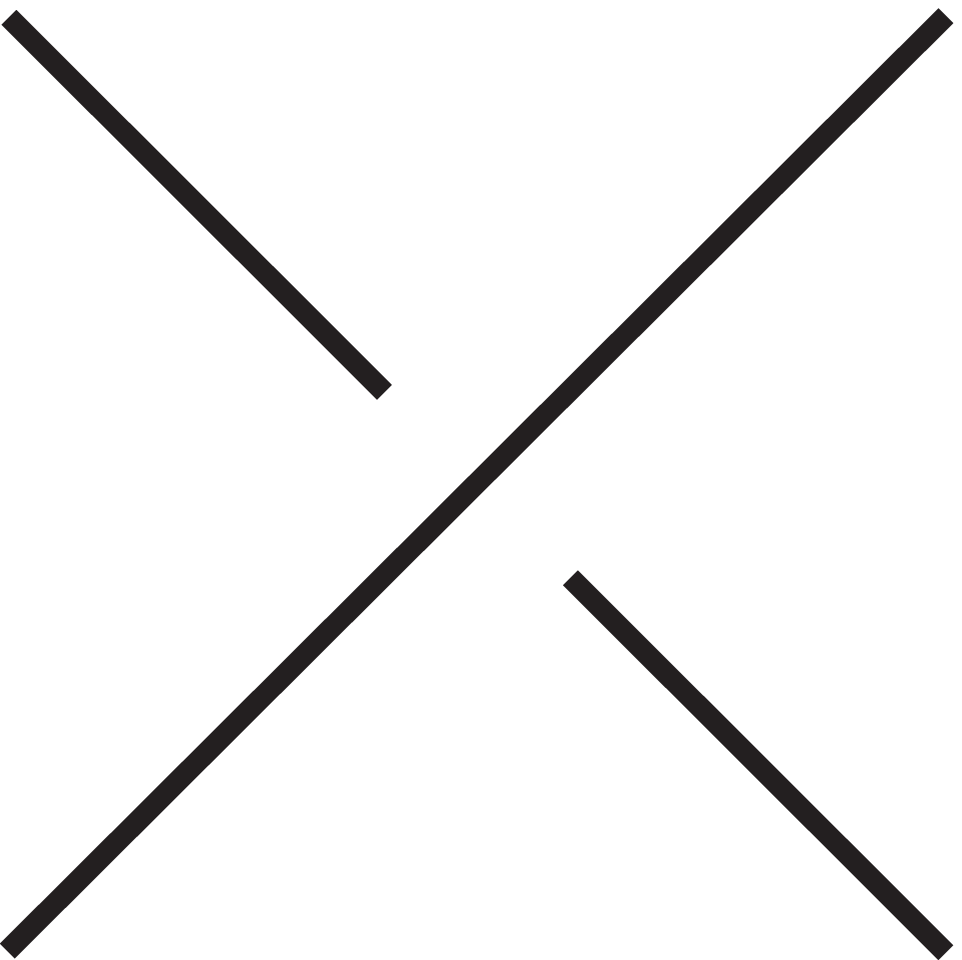}}\to
 \raisebox{-0.25\height}{\includegraphics[width=0.5cm]{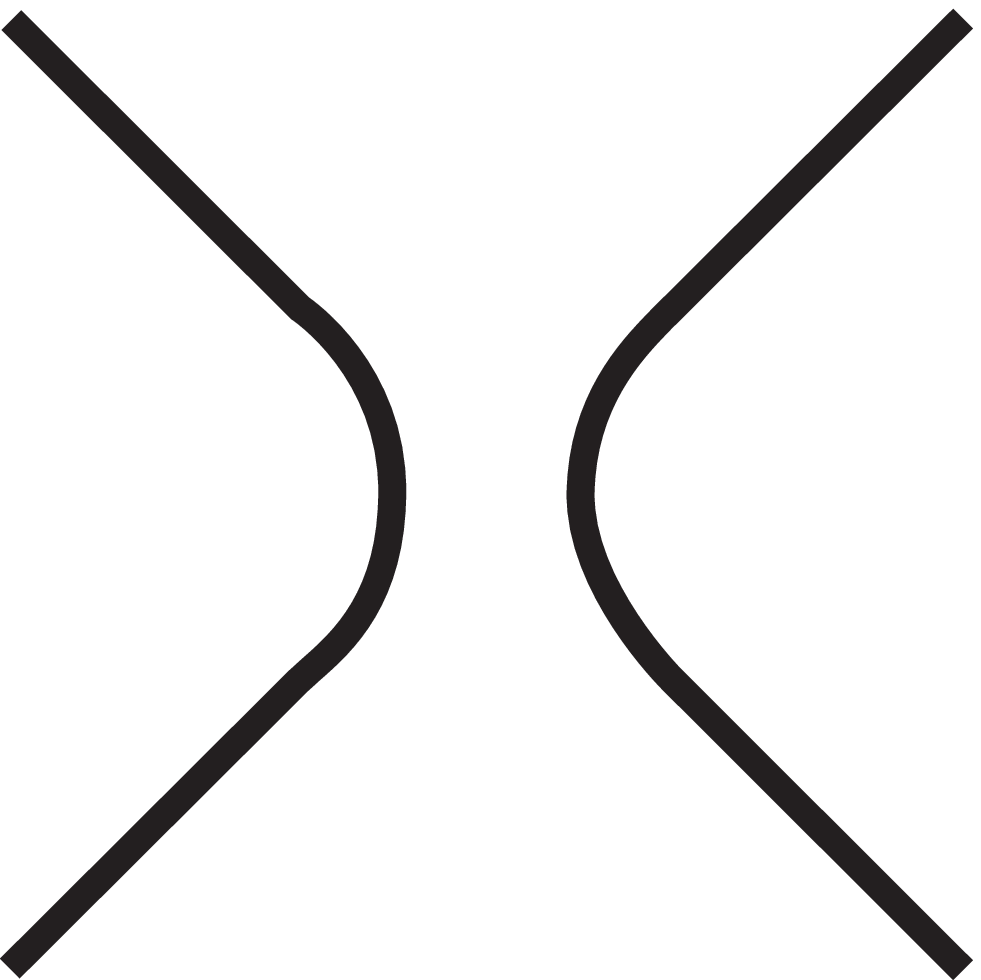}}$,
 $B\colon\raisebox{-0.25\height}{\includegraphics[width=0.5cm]{skcrossr.eps}}\to
 \raisebox{-0.25\height}{\includegraphics[width=0.5cm]{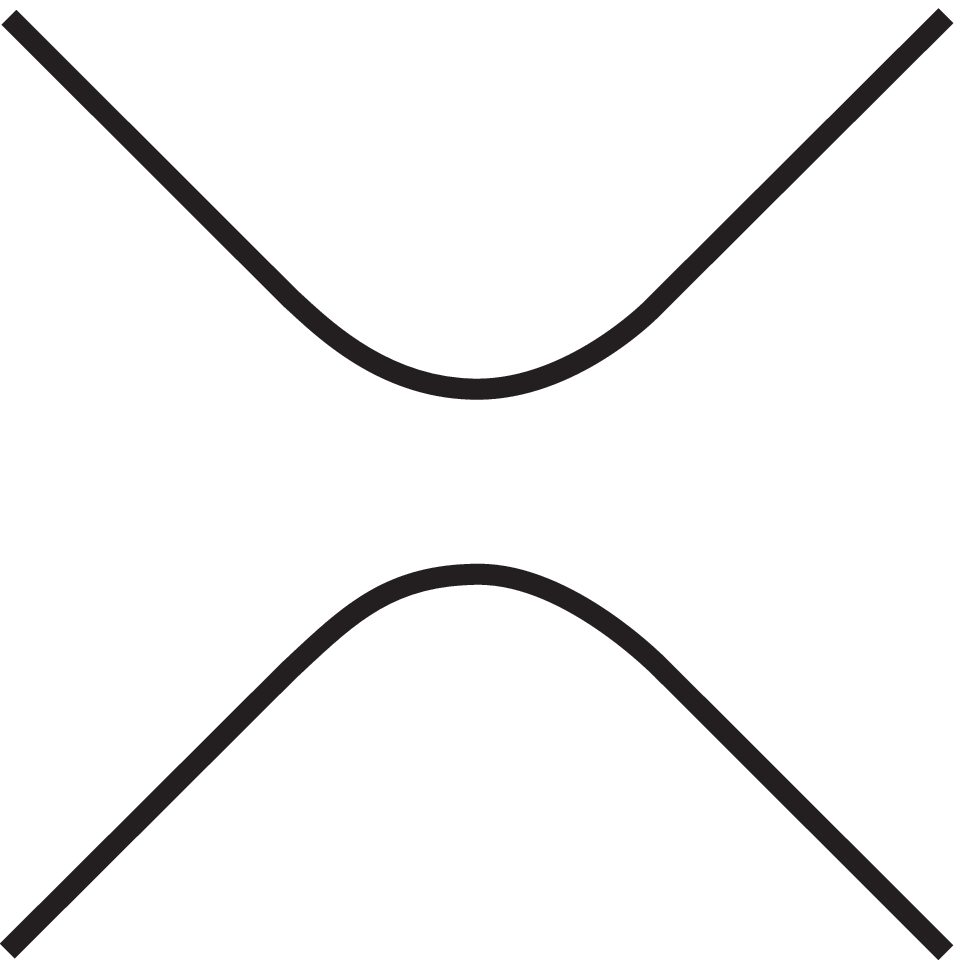}}$.

 \begin {figure} \centering\includegraphics[width=300pt]{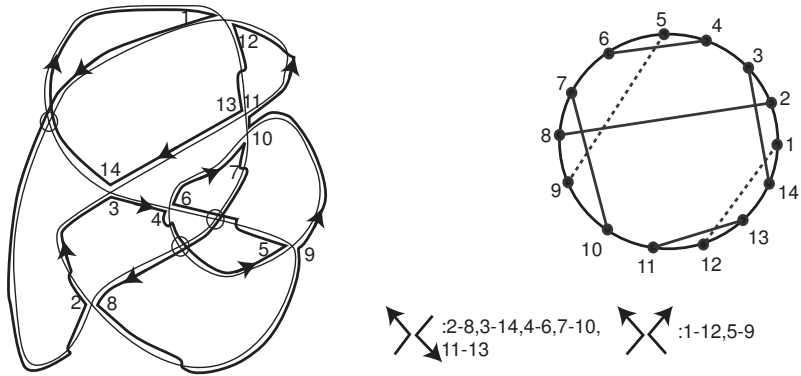}
\caption{Rotating Circuit Shown by a Thick Line; Chord
Diagram}\label{rotations}
 \end {figure}

Conversely, having a connected virtual diagram $K$, one can get a
labeled chord diagram $D_C(K)$, see Fig.~\ref{rotations}. Indeed,
one takes a {\em circuit} $C$ of $K$ which is a map from $S^{1}$ to
the projection of $K$. This map is bijective outside classical and
virtual crossings, has exactly two preimages at each classical and
virtual crossing, goes transversally at each virtual crossing and
turns from an half-edge to an adjacent (non-opposite) half-edge at
each classical crossing. Connecting the two preimages of a classical
crossing by a chord we get a chord diagram, where the sign of the
chord is $+$ if the circuit locally agrees with the $A$-smoothing,
and $-$ if it agrees with the $B$-smoothing, and the framing of a
chord is $0$ (resp., $1$) if two opposite half-edges have the
opposite (resp., the same) orientation. It can be easily checked
that this operation is indeed inverse to the operation of
constructing a virtual link out of a chord diagram: if we take a
chord diagram $D$, and construct a virtual diagram $K(D)$ out of it,
then for some circuit $C$ the chord diagram $D_C(K(D))$ will
coincide with $D$. The rule for setting classical crossings here
agrees with the rule described above. This proves the following

  \begin {theorem}[\cite {Ma3}]
For any connected virtual diagram $L$ there is a certain labeled
chord diagram $D$ such that $L=K(D)$.
  \end {theorem}

The Reidemeister moves on virtual diagrams generate the Reidemeister
moves on labeled chord diagrams~\cite {IM1,IM2}.

\subsection{Reidemeister Moves for Looped Interlacement Graphs and Graph-links}

Now we are describing moves on graphs obtained from virtual diagrams
by using rotating circuit~\cite {IM1,IM2} and the Gauss
circuit~\cite {TZ}. These moves in both cases will correspond to the
``real'' Reidemeister moves on diagrams. Then we shall extend these
moves to all graphs (not only to realisable ones). As a result we
get new objects, a {\em graph-link} and a {\em homotopy class of
looped interlacement graphs}, in a way similar to the generalisation
of classical knots to virtual knots: the passage from realisable
Gauss diagrams (classical knots) to arbitrary chord diagrams leads
to the concept of a virtual knot, and the passage from realisable
(by means of chord diagrams) graphs to arbitrary graphs leads to the
concept of two new objects: graph-links and homotopy classes of
looped interlacement graphs (here `looped' corresponds to the
writher number, if the writher number is -1 then the corresponding
vertex has a loop). To construct the first object we shall use
simple labeled graphs, and for the second one we shall use
(unlabeled) graphs without multiple edges, but loops are allowed.

 \begin {definition}
A graph is {\em labeled} if every vertex $v$ of it is endowed with a
pair $(a,\alpha)$,  where $a\in\{0,1\}$ is the framing of $v$, and
$\alpha\in\{\pm\}$ is the sign of $v$. Let $D$ be a labeled chord
diagram $D$. The {\em labeled intersection graph}, cf.~\cite{CDL},
$G(D)$ of $D$ is the labeled graph: 1) whose vertices are in
one-to-one correspondence with chords of $D$, 2) the label of each
vertex corresponding to a chord coincides with that of the chord,
and 3) two vertices are connected by an edge if and only if the
corresponding chords are linked.
 \end {definition}

 \begin {definition}
A simple (labeled) graph $H$ is called {\em realisable} if there is
a (labeled) chord diagram $D$ such that $H=G(D)$.
 \end {definition}

The following lemma is evident.

 \begin {lemma}\label {lem:co_re}
A simple {\em(}labeled\/{\em)} graph is realisable if and only if
each its connected component is realisable.
 \end {lemma}

 \begin {definition}
Let $G$ be a graph and let $v\in V(G)$. The set of all vertices
adjacent to $v$ is called the {\em neighbourhood of a vertex} $v$
and denoted by $N(v)$ or $N_G(v)$.
 \end {definition}

Let us define two operations on simple unlabeled graphs.

 \begin {definition} (Local Complementation)
Let $G$ be a graph. The {\em local complementation} of $G$ at $v\in
V (G)$ is the operation which toggles adjacencies between $a,b\in
N(v)$, $a\ne b$, and doesn't change the rest of $G$. Denote the
graph obtained from $G$ by the local complementation at a vertex $v$
by $\operatorname{LC}(G;v)$.
 \end {definition}

 \begin {definition} (Pivot)
Let $G$ be a graph with distinct vertices $u$ and $v$. The {\em
pivoting operation} of a graph $G$ at $u$ and $v$ is the operation
which toggles adjacencies between $x,\,y$ such that
$x,\,y\notin\{u,v\}$, $x\in N(u),\,y\in N(v)$ and either $x\notin
N(v)$ or $y\notin N(u)$, and doesn't change the rest of $G$. Denote
the graph obtained from $G$ by the pivoting operation at vertices
$u$ and $v$ by $\operatorname{piv}(G;u,v)$.
 \end {definition}

\begin{figure} \centering\includegraphics[width=200pt]{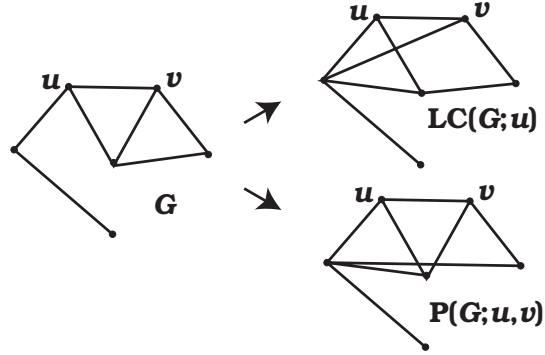}
  \caption{Local Complementation and Pivot} \label{lcp}
 \end{figure}

 \begin {example}
In Fig.~\ref {lcp} the graphs $G$, $\operatorname{LC}(G;u)$ and
$\operatorname{piv}(G;u,v)$ are depicted.
 \end {example}

The following lemma can be easily checked.

 \begin {lemma}
If $u$ and $v$ are adjacent then there is an isomorphism
 $$
\operatorname{piv}(G;u,v)\cong\operatorname{LC}(\operatorname{LC}(\operatorname{LC}(G;u);v);u).
 $$
 \end {lemma}

Let us define graph-moves by considering intersection graphs of
chord diagrams constructed by using {\em a rotating circuit}, and
these moves correspond to the Reidemeister moves on virtual
diagrams. As a result we obtain a new object --- an equivalence
class of labeled graphs under formal moves. These moves were defined
in~\cite {IM1,IM2}.

 \begin {definition}
$\Omega_g 1$. The first Reidemeister graph-move is an
addition/removal of an isolated vertex labeled $(0,\alpha)$,
$\alpha\in\{\pm\}$.

$\Omega_g 2$. The second Reidemeister graph-move is an
addition/removal of two non-adjacent (resp., adjacent) vertices
having $(0,\pm\alpha)$ (resp., $(1,\pm\alpha)$) and the same
adjacencies with other vertices.

$\Omega_g 3$. The third Reidemeister graph-move is defined as
follows. Let $u,\,v,\,w$ be three vertices of $G$ all having label
$(0,-)$ so that $u$ is adjacent only to $v$ and $w$ in $G$. Then we
only change the adjacency of $u$ with the vertices $t\in
N(v)\setminus N(w)\bigcup N(w)\setminus N(v)$ (for other pairs of
vertices we do not change their adjacency). In addition, we switch
the signs of $v$ and $w$ to $+$. The inverse operation is also
called the third Reidemeister graph-move.

$\Omega_g 4$. The fourth graph-move for $G$ is defined as follows.
We take two adjacent vertices $u$ and $v$ labeled $(0,\alpha)$ and
$(0,\beta)$ respectively. Replace $G$ with
$\operatorname{piv}(G;u,v)$ and change signs of $u$ and $v$ so that
the sign of $u$ becomes $-\beta$ and the sign of $v$ becomes
$-\alpha$.

$\Omega_g 4'$. In this fourth graph-move we take a vertex $v$ with
the label $(1,\alpha)$. Replace $G$ with $\operatorname{LC}(G;v)$
and change the sign of $v$ and the framing for each $u\in N(v)$.
 \end {definition}

The comparison of the graph-moves with the Reidemeister moves yields
the following theorem.

 \begin {theorem}\label {th:link_graph}
Let $K_1$ and $K_2$ be two connected virtual diagrams, and let $G_1$
and $G_2$ be two labeled intersection graphs obtained from $K_1$ and
$K_2$, respectively. If $K_1$ and $K_2$ are equivalent in the class
of connected diagrams then $G_1$ and $G_2$ are obtained from one
another by a sequence of $\Omega_g 1 - \Omega_g 4'$ graph-moves.
 \end {theorem}

  \begin {definition}
A {\em graph-link} is an equivalence class of simple labeled graphs
modulo $\Omega_g 1 - \Omega_g 4'$ graph-moves.
 \end {definition}

Let $D_G(K)$ be the Gauss diagram of a virtual diagram $K$. Let us
construct the graph obtained from the intersection graph of $D_G(K)$
by adding loops to vertices corresponding to chords with negative
writhe number~\cite {TZ}. We refer to this graph as the {\em looped
interlacement graph} or the {\em looped graph}. Let us construct the
moves on graphs. These moves are similar to the moves for
graph-links and also correspond to the Reidemeister moves on virtual
diagrams.

 \begin {definition}\label {def:mov_ga}
$\Omega 1$. The first Reidemeister move for looped interlacement
graphs is an addition/removal of an isolated looped or unlooped
vertex.

$\Omega 2$. The second Reidemeister move for looped interlacement
graphs is an addition/removal of two vertices having the same
adjacencies with other vertices and, moreover, one of which is
looped and the other one is unlooped.

$\Omega 3$. The third Reidemeister move for looped interlacement
graphs is defined as follows. Let $u,\,v,\,w$ be three vertices such
that $v$ is looped, $w$ is unlooped, $v$ and $w$ are adjacent, $u$
is adjacent to neither $v$ nor $w$, and every vertex
$x\notin\{u,\,v,\,w\}$ is adjacent to either $0$ or precisely two of
$u,\,v,\,w$. Then we only remove all three edges $uv$, $uw$ and
$vw$.  The inverse operation is also called the third Reidemeister
move.
 \end {definition}

\begin{figure}
\centering\includegraphics[width=300pt]{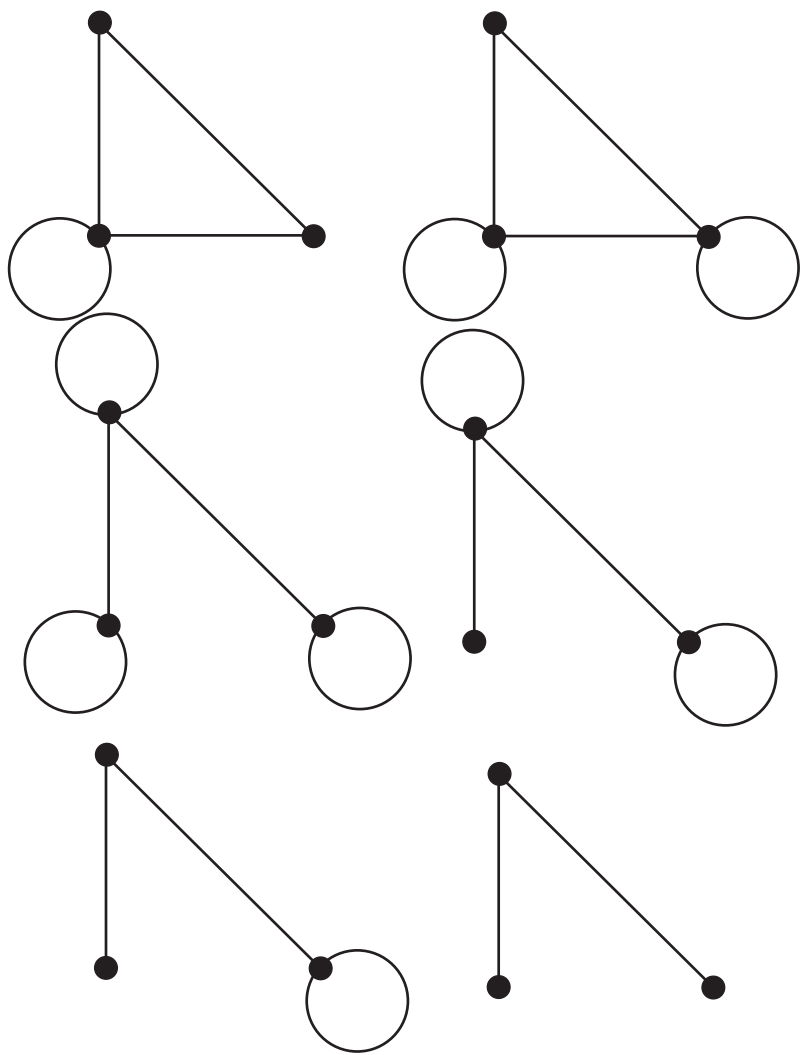} \caption{The
Possible Configurations of the 3-rd Reidemeister move} \label{reid}
 \end{figure}

 \begin {remark}\label {rk:loop3move}
The two third Reidemeister moves do not exhaust all the
possibilities for representing the third Reidemeister move on Gauss
diagrams, see~\cite {TZ}. It can be shown that all the other
versions of the third Reidemeister move, which involve toggling the
non-loop edges in one of the six pictured configurations in
Figure~\ref {reid} (every vertex outside the picture must have
either $0$ or precisely $2$ neighbours among the three vertices that
are pictured), are combinations of the second and third Reidemeister
moves, see~\cite {Ost} for details.
 \end {remark}

 \begin {definition}
We call an equivalence class of graphs (without multiple edges, but
loops are allowed) modulo the three moves listed in Definition~2.11
a {\em homotopy class} of looped interlacement graphs.
 \end {definition}

 \begin {remark}
Looped interlacement graphs encode only knot diagrams but
graph-links can encode virtual diagrams with any number of
components. The approach using a rotating circuit has an advantage
in this sense. In~\cite {Tr1} L.~Traldi introduced the notion of a
marked graph by considering any Euler tour (we have vertices which
we go transversally and in which we rotate). In the paper we don't
consider ``mixes circuits''.
 \end {remark}

\section{Graph-Links and Homotopy Classes of Looped Graphs: Their Equivalence}

The main goal of this section is to construct an equivalence between
the set of graph-knots, see ahead, and the set of homotopy classes
of looped graphs such that the graph-knot and the homotopy class of
looped graphs constructed from a given knot are related by the
equivalence.

Before constructing an equivalence we need some definitions and
statements.

 \begin {definition}
Let $G$ be a labeled graph on vertices from the enumerated set
$V(G)=\{v_1,\dots,v_n\}$. Define the {\em adjacency matrix}
$A(G)=(a_{ij})$ of $G$ over $\Z_2$ as follows: $a_{ii}$ is equal to
the framing of $v_i$, $a_{ij}=1$, $i\ne j$, if and only if $v_i$ is
adjacent to $v_j$ and $a_{ij}=0$ otherwise.
 \end {definition}

 \begin {remark}
Throughout the paper all the matrices are over $\mathbb{Z}_2$.
Therefore, $\corank$'s and $\det$'s are calculated over
$\mathbb{Z}_2$.
 \end {remark}

If $G$ and $G'$ represent the same graph-link then
$\operatorname{corank}_{\mathbb{Z}_2}(A(G)+E)=\operatorname{corank}_{\mathbb{Z}_2}(A(G')+E)$,
where $E$ is identity matrix,  see~\cite {IM1, IM2}.

 \begin {definition}\label {def:num_com}
Let us define the {\em number of components} in a graph-link
$\mathfrak{F}$ as $\operatorname{corank}_{\mathbb{Z}_2}(A(G)+E)+1$,
here $G$ is a representative of $\mathfrak{F}$. A graph-link
$\mathfrak{F}$ with $\operatorname{corank}_{\mathbb{Z}_2}(A(G)+E)=0$
for any representative $G$ of $\mathfrak{F}$ is called a {\em
graph-knot}.
 \end {definition}

Let $\operatorname{corank}_{\mathbb{Z}_2}(A(G)+E)=0$,
$B_i(G)=A(G)+E+E_{ii}$ (all elements of $E_{ii}$ except for the one
in the $i$-th column and $i$-th row which is one are $0$) for each
vertex $v_i\in V(G)$.

 \begin {definition}\label {def:wr_num}
The {\em writhe number} $w_i$ of $G$ (with
$\operatorname{corank}_{\mathbb{Z}_2}(A(G)+E)=0$) at $v_i$ is
$w_i=(-1)^{\operatorname{corank}_{\mathbb{Z}_2}
B_i(G)}\operatorname{sign}v_i$ and the {\em writhe number} of $G$ is
  $$
w(G)=\sum\limits_{i=1}^nw_i.
  $$
 \end {definition}

If $G$ is a realisable graph by a chord diagram and, therefore, by a
virtual diagram then $w_i$ is the ``real'' writhe number of the
crossing corresponding to $v_i$.

By $\widehat{C}_{i,j,\ldots,k}$ denote the matrix obtained from a
matrix $C$ by deleting the $i,j,\ldots,k$-th rows and the
$i,j,\ldots,k$-th columns. We shall write
$\widehat{B}_{i,j,\ldots,k}(G)$ instead of
$\widehat{B(G)}_{i,j,\ldots,k}$.

 \begin {lemma}
$w_i(G)=(-1)^{\operatorname{corank}_{\mathbb{Z}_2}
\widehat{B}_i(G)+1}\operatorname{sign}v_i$.
 \end {lemma}

 \begin {proof}
Without loss of generality we prove the statement of the lemma for
$w_1$. Let
 $$
A(G)+E=\left(\begin {array}{cc} a & \bb^\top \\
\bb & C
 \end
 {array}\right)
 $$
and
 $$
\corank_{\mathbb{Z}_2}(A(G)+E)=0\Longleftrightarrow\det(A(G)+E)=1,
 $$
where bold characters indicate column vectors. We have
 $$
\det B_1(G)=\det\left(\begin {array}{cc} a+1 & \bb^\top \\
\bb & C
 \end
 {array}\right)=\det\left(\begin {array}{cc} a & \bb^\top \\
\bb & C
 \end
 {array}\right)+\det\left(\begin {array}{cc} 1 & \0^\top \\
\bb & C
 \end
 {array}\right)
 $$
and
 $$
\widehat{B}_1(G)=C,\quad \det(B_1(G))=\det A(G)+\det
C=1+\det\widehat{B}_1(G).
 $$
The last equality gives us the statement of the lemma.
 \end {proof}

 \begin {lemma}\label {lem:wrfr_si}
Let $G$ be labeled graph with $\det B(G)=1$ and
$B(G)^{-1}=(b^{ij})$. Then $b^{ii}=\dfrac{1-w_i(G)\sign v_i}2$.
 \end {lemma}

 \begin {proof}
We have
 $$
w_i(G)=(-1)^{\corank_{\mathbb{Z}_2}\widehat{B}_i(G)+1}\sign
v_i\Longleftrightarrow w_i(G)\sign
v_i=(-1)^{\corank_{\mathbb{Z}_2}\widehat{B}_i(G)+1}
 $$
 $$
\Longleftrightarrow\corank_{\mathbb{Z}_2}\widehat{B}_i(G)=\dfrac{w_i(G)\sign
v_i+1}2\Longleftrightarrow b^{ii}=\det\widehat{B}_i(G)
 $$
 $$
=1-\corank_{\mathbb{Z}_2}\widehat{B}_i(G)=\dfrac{1-w_i(G)\sign
v_i}2.
 $$
 \end {proof}

 \begin {definition}
Define the {\em adjacency matrix} $A(L)=(a_{ij})$ over $\Z_2$ for a
looped graph $L$ with enumerated set of vertices as: $a_{ii}=1$ if
and only if the vertex numbered $i$ is looped and $a_{ii}=0$
otherwise; $a_{ij}=1$, $i\ne j$, if and only if the vertex with the
number $i$ is adjacent to the vertex with the number $j$ and
$a_{ij}=0$ otherwise.
 \end {definition}

 \begin {definition}
We say that an $n\times n$ matrix $A=(a_{ij})$ {\em coincides with
an $n\times n$ matrix $B=(b_{kl})$ up to diagonal elements} if
$a_{ij}=b_{ij}$, $i\ne j$.
 \end {definition}

  \begin {lemma}[\cite {Ily}]\label {lem:nondeg}
Let $A$ be a symmetric matrix over $\mathbb{Z}_2$. Then there exists
a symmetric matrix $\widetilde{A}$ with $\det\widetilde{A}=1$ equal
to $A$ up to diagonal elements.
 \end {lemma}

The main theorem of the paper is the following one.

 \begin {theorem}\label {th:gk_hl}
There is a one-to-one correspondence between the set of all
graph-knots and the set of all homotopy classes of looped graphs.
Moreover, if $K$ is a knot, and $\mathfrak{F}$ is a graph-knot
constructed from $K$ and $\mathfrak{L}$ is the homotopy class of
looped graphs constructed from $K$ then $\mathfrak{F}$ and
$\mathfrak{L}$ are related by this correspondence.
 \end {theorem}

To prove this theorem we construct a map from the set of all
graph-knots to the set of all homotopy classes of looped graphs. We
shall show that this map has the inverse map.

Let us construct the map $\chi$ from the set of all graph-knots to
the set of all homotopy classes of looped graphs. Let $\mathfrak{F}$
be a graph-knot and let $G$ be its representative. Let us consider
the simple graph $H$ having the adjacency matrix coinciding with
$(A(G)+E)^{-1}$ up to diagonal elements and construct the graph
$L(G)$ from $H$ by just adding loops to any vertex of $H$
corresponding to a vertex of $G$ with the negative writhe number. By
definition, put $\chi(\mathfrak{F})=\mathfrak{L}$, here
$\mathfrak{L}$ is the homotopy class of $L(G)$.

 \begin {theorem}\label {th:grkn_ho}
The map $\chi$ from the set of all graph-knots to the set of all
homotopy classes of looped graphs defined by
$\chi(\mathfrak{F})=\mathfrak{L}$, here $\mathfrak{L}$ is the
homotopy class of $L(G)$ and $G$ is a representative of
$\mathfrak{F}$, is well-defined.
 \end {theorem}

 \begin {proof}
Let $G_1,\,G_2$ be two representatives of $\mathfrak{F}$, and let
$B(G_i)=A(G_i)+E$, $B(G_i)^{-1}=(b^{kl}_i)$. We have to show that
the homotopy classes of  $L(G_1)$ and $L(G_2)$ are the same, i.e.\
$L_1=L(G_1)$ and $L_2=L(G_2)$ are related to each other by the
Reidemeister moves on looped graphs.

Let us consider four cases.

$1)$ We know that if $G_1$ and $G_2$ are obtained from each other by
$\Omega_g4$ and/or $\Omega_g4'$ (these moves correspond to changing
of a rotating circuit in the case of realisable graphs) then the
writhe numbers of corresponding vertices of $G_1$ and $G_2$ are the
same, see~~\cite {IM1,IM2}, and the graphs $L_1$ and $L_2$ are
isomorphic up to loops, see~\cite {Ily}. Therefore, $L_1$ and $L_2$
are isomorphic.

$2)$ Let $G_1$ and $G_2$ be obtained from each other by $\Omega_g1$
(we remove the vertex with the number $1$). We have
 $$
B(G_1)=\left(\begin {array} {cc}
 1 & {\mathbf 0}^\top\\
 {\mathbf 0} & A(G_2)+E
 \end {array}\right),\qquad B(G_2)=A(G_2)+E,
 $$
 $$
B(G_1)^{-1}=\left(\begin {array} {cc}
 1 & {\mathbf 0}^\top\\
 {\mathbf 0} & (A(G_2)+E)^{-1}
 \end {array}\right),\qquad B(G_2)^{-1}=(A(G_2)+E)^{-1},
 $$
where $\mathbf 0$ indicates column vector with all entries $0$.
Therefore, $L_2$ is obtained from $L_1$ by $\Omega1$.

$3)$ Let $G_1$ and $G_2$ be obtained from each other by $\Omega_g2$
(we remove the vertex with the numbers $1$ and $2$). We have two
cases. The first case
 $$
B(G_1)=\left(\begin {array} {ccc}
 1 & 0 & {\mathbf a}^\top\\
 0 & 1 & {\mathbf a}^\top\\
 {\mathbf a} & {\mathbf a}& A(G_2)+E
 \end {array}\right),\qquad B(G_2)=A(G_2)+E,
 $$
and the second case
 $$
B(G_1)=\left(\begin {array} {ccc}
 0 & 1 & {\mathbf a}^\top\\
 1 & 0 & {\mathbf a}^\top\\
 {\mathbf a} & {\mathbf a}& A(G_2)+E
 \end {array}\right),\qquad B(G_2)=A(G_2)+E.
 $$
Let us consider only the first case. We know that
$w_1(G_1)=-w_2(G_1)$ in $G_1$~\cite {IM1,IM2}. Moreover, as
 $$
\det\left(\begin {array} {ccc}
 1 & 0 & \widetilde{{\mathbf a}}^\top\\
 0 & 1 & \widetilde{{\mathbf a}}^\top\\
 \widetilde{{\mathbf a}} & \widetilde{{\mathbf a}} & C
 \end {array}\right)=\det\left(\begin {array} {ccc}
 1 & 0 & \widetilde{{\mathbf a}}^\top\\
 1 & 1 & {\mathbf 0}^\top\\
 \widetilde{{\mathbf a}} & \widetilde{{\mathbf a}} & C
 \end {array}\right)=\det\left(\begin {array} {ccc}
 1 & 0 & \widetilde{{\mathbf a}}^\top\\
 1 & 1 & {\mathbf 0}^\top\\
 {\mathbf 0} & {\mathbf 0}& C
 \end {array}\right)=\det C,
 $$
then
 $$
B(G_1)^{-1}=\left(\begin {array} {ccc}
 b & c & {\mathbf d}^\top\\
 c & b & {\mathbf d}^\top\\
 {\mathbf d} & {\mathbf d}& (A(G_2)+E)^{-1}
 \end {array}\right),\quad B(G_2)^{-1}=(A(G_2)+E)^{-1}.
 $$
This means that $L_1$ and $L_2$ are obtained from each other by
$\Omega 2$.

$4)$ Now assume that $G_1$ and $G_2$ are obtained from each other by
$\Omega_g3$. The corresponding vertices of $G_1$ and $G_2$ under
$\Omega_g3$ have the same numbers (in~\cite {IM1} we have another
enumeration). We shall prove that $L_1$ and $L_2$ are obtained from
each other by a sequence of $\Omega 2$, $\Omega 3$ moves.

We have
 $$
B(G_1)=\left(\begin {array} {cccc}
 1 & 1 & 1 & {\mathbf 0}^\top\\
 1 & 1 & 0 & {\mathbf a}^\top\\
 1 & 0 & 1 & {\mathbf b}^\top\\
 {\mathbf 0} & {\mathbf a} & {\mathbf b} & C
 \end {array}\right),
 $$
 $$
1=\det B(G_1)=\det\left(\begin {array} {cccc}
 1 & 1 & 1 & {\mathbf 0}^\top\\
 1 & 1 & 0 & {\mathbf a}^\top\\
 1 & 0 & 1 & {\mathbf b}^\top\\
 {\mathbf 0} & {\mathbf a} & {\mathbf b} & C
 \end {array}\right)=\det\left(\begin {array} {ccc}
 0 & 1 & {\mathbf a}^\top\\
 1 & 0 & {\mathbf b}^\top\\
 {\mathbf a} & {\mathbf b} & C
 \end {array}\right),
 $$
 $$
B(G_2)=\left(\begin {array} {cccc}
 1 & 0 & 0 & ({\mathbf a}+{\mathbf b})^\top\\
 0 & 1 & 0 & {\mathbf b}^\top\\
 0 & 0 & 1 & {\mathbf a}^\top\\
 {\mathbf a}+{\mathbf b} & {\mathbf b} & {\mathbf a} & C
 \end {array}\right),\qquad \det B(G_2)=1
 $$

Let us show that we have a structure either for $v_1,\,v_2,\,v_3\in
V(L_1)$ or $\widetilde{v}_1,\,\widetilde{v}_2,\,\widetilde{v}_3\in
V(L_2)$ as in Figure~\ref {reid}.

We have (\cite {IM1})
 $$
w_i(G_1)=w_i(G_2),\quad i=1,\,2,\,3,
 $$
 $$
\det\widehat{B}_1(G_1)=\det\left(\begin {array} {ccc}
 1 & 0 & {\mathbf a}^\top\\
 0 & 1 & {\mathbf b}^\top\\
{\mathbf a} & {\mathbf b} & C
 \end {array}\right)
 $$
 $$
=\det\left(\begin {array} {ccc}
 1 & 0 & {\mathbf a}^\top\\
 1 & 1 & {\mathbf b}^\top\\
{\mathbf 0} & {\mathbf b} & C
 \end {array}\right)+\det\left(\begin {array} {ccc}
 0 & 0 & {\mathbf a}^\top\\
 1 & 1 & {\mathbf b}^\top\\
{\mathbf a} & {\mathbf b} & C
 \end {array}\right)
 $$
 $$
=\det\left(\begin {array} {ccc}
 1 & 0 & {\mathbf a}^\top\\
 1 & 1 & {\mathbf b}^\top\\
{\mathbf 0} & {\mathbf b} & C
 \end {array}\right)+\det\left(\begin {array} {ccc}
 0 & 1 & {\mathbf a}^\top\\
 1 & 1 & {\mathbf b}^\top\\
{\mathbf a} & {\mathbf 0} & C
 \end {array}\right)+\det\left(\begin {array} {ccc}
 0 & 1 & {\mathbf a}^\top\\
 1 & 0 & {\mathbf b}^\top\\
{\mathbf a} & {\mathbf b} & C
 \end {array}\right),
 $$
 $$
\det\widehat{B}_2(G_1)=\det\left(\begin {array} {ccc}
 1 & 1 & {\mathbf 0}^\top\\
 1 & 1 & {\mathbf b}^\top\\
 {\mathbf 0} & {\mathbf b} & C
 \end {array}\right)=\det\left(\begin {array} {cc}
 0 & {\mathbf b}^\top\\
 {\mathbf b} & C
 \end {array}\right),
 $$
 $$
\det\widehat{B}_3(G_1)=\det\left(\begin {array} {ccc}
 1 & 1 & {\mathbf 0}^\top\\
 1 & 1 & {\mathbf a}^\top\\
 {\mathbf 0} & {\mathbf a} & C
 \end {array}\right)=\det\left(\begin {array} {cc}
 0 & {\mathbf a}^\top\\
 {\mathbf a} & C
 \end {array}\right),
 $$
 $$
b_1^{12}=\det\left(\begin {array} {ccc}
 1 & 0 & {\mathbf a}^\top\\
 1 & 1 & {\mathbf b}^\top\\
 {\mathbf 0} & {\mathbf b} & C
 \end {array}\right)=\det\left(\begin {array} {cc}
  1 & {\mathbf a}^\top+{\mathbf b}^\top\\
 {\mathbf b} & C
 \end {array}\right)
 $$
 $$
=\det\left(\begin {array} {cc}
  1 & {\mathbf a}^\top\\
 {\mathbf b} & C
 \end {array}\right)+\det\left(\begin {array} {cc}
 0 & {\mathbf b}^\top\\
 {\mathbf b} & C
 \end {array}\right),
 $$
 $$
b_1^{13}=\det\left(\begin {array} {ccc}
 1 & 0 & {\mathbf a}^\top\\
 1 & 1 & {\mathbf b}^\top\\
 {\mathbf 0} & {\mathbf a} & C
 \end {array}\right)=\det\left(\begin {array} {cc}
  1 & {\mathbf a}^\top+{\mathbf b}^\top\\
 {\mathbf a} & C
 \end {array}\right)
 $$
 $$
=\det\left(\begin {array} {cc}
  1 & {\mathbf b}^\top\\
 {\mathbf a} & C
 \end {array}\right)+\det\left(\begin {array} {cc}
 0 & {\mathbf a}^\top\\
 {\mathbf a} & C
 \end {array}\right),
 $$
 $$
b_1^{23}=\det\left(\begin {array} {ccc}
 1 & 1 & {\mathbf 0}^\top\\
 1 & 0 & {\mathbf b}^\top\\
 {\mathbf 0} & {\mathbf a} & C
 \end {array}\right)=\det\left(\begin {array} {cc}
 1 & {\mathbf b}^\top\\
 {\mathbf a} & C
 \end {array}\right),
 $$
 $$
b_2^{12}=\det\left(\begin {array} {ccc}
 0 & 0 & {\mathbf b}^\top\\
 0 & 1 & {\mathbf a}^\top\\
 {\mathbf a}+{\mathbf b} & {\mathbf a} & C
 \end {array}\right)=\det\left(\begin {array} {ccc}
 0 & 0 & {\mathbf b}^\top\\
 1 & 1 & {\mathbf a}^\top\\
 {\mathbf b} & {\mathbf a} & C
 \end {array}\right)
 $$
 $$
=\det\left(\begin {array} {ccc}
 0 & 1 & {\mathbf b}^\top\\
 1 & 0 & {\mathbf a}^\top\\
 {\mathbf b} & {\mathbf a} & C
 \end {array}\right)+\det\left(\begin {array} {ccc}
 0 & 1 & {\mathbf b}^\top\\
 1 & 1 & {\mathbf a}^\top\\
 {\mathbf b} & {\mathbf 0} & C
 \end {array}\right)=1+b^{12}_1,
 $$
 $$
b_2^{13}=\det\left(\begin {array} {ccc}
 0 & 1 & {\mathbf b}^\top\\
 0 & 0 & {\mathbf a}^\top\\
 {\mathbf a}+{\mathbf b} & {\mathbf b} & C
 \end {array}\right)=\det\left(\begin {array} {ccc}
 1 & 1 & {\mathbf b}^\top\\
 0 & 0 & {\mathbf a}^\top\\
 {\mathbf a} & {\mathbf b} & C
 \end {array}\right)
 $$
 $$
=\det\left(\begin {array} {ccc}
 1 & 0 & {\mathbf b}^\top\\
 0 & 1 & {\mathbf a}^\top\\
 {\mathbf a} & {\mathbf b} & C
 \end {array}\right)+\det\left(\begin {array} {ccc}
 1 & 1 & {\mathbf b}^\top\\
 0 & 1 & {\mathbf a}^\top\\
 {\mathbf a} & {\mathbf 0} & C
 \end {array}\right)=1+b^{13}_1,
 $$
 $$
b_2^{23}=\det\left(\begin {array} {ccc}
 1 & 0 & ({\mathbf a}+{\mathbf b})^\top\\
 0 & 0 & {\mathbf a}^\top\\
 {\mathbf a}+{\mathbf b} & {\mathbf b} & C
 \end {array}\right)=\det\left(\begin {array} {ccc}
 1 & 0 & {\mathbf b}^\top\\
 0 & 0 & {\mathbf a}^\top\\
 {\mathbf a} & {\mathbf b} & C
 \end {array}\right)
 $$
 $$
=\det\left(\begin {array} {ccc}
 1 & 0 & {\mathbf b}^\top\\
 0 & 1 & {\mathbf a}^\top\\
 {\mathbf a} & {\mathbf b} & C
 \end {array}\right)+\det\left(\begin {array} {ccc}
 1 & 0 & {\mathbf b}^\top\\
 0 & 1 & {\mathbf a}^\top\\
 {\mathbf a} & {\mathbf 0} & C
 \end {array}\right)=1+b_2^{23}.
 $$
 $$
b^{12}_1=b^{23}_1+\det\widehat{B}_2(G_1),\quad
b^{13}_1=b^{23}_1+\det\widehat{B}_3(G_1),
 $$
 $$
b^{12}_1+b^{13}_1=\det\widehat{B}_1(G_1)+1,
 $$
 $$
b^{12}_2=b^{12}_1+1,\qquad b^{13}_2=b^{13}_1+1,\qquad
b^{23}_2=b^{23}_1+1.
$$
It is not difficult to show that the last equalities guarantee us
that either $v_1,\,v_2,\,v_3\in V(L_1)$ or
$\widetilde{v}_1,\,\widetilde{v}_2,\,\widetilde{v}_3\in V(L_2)$ have
a structure as in Figure~\ref {reid}. The structure of the other
triple is obtained from the first triple by toggling the non-loop
edges.

Denote by $\widehat{{\mathbf f}}^i$ the column vector obtained from
${\mathbf f}$ by deleting the $i$-th element and denote by
$\widehat{C}_j$ (resp., $\widehat{C}_j^i$) the matrix obtained from
$C$ by deleting the $j$-th column (resp., the $i$-th row and the
$j$-th column).

Other equalities for $i,j>3$ are as follows:
 $$
b^{ij}_1=\det\left(\begin {array} {cccc}
 1 & 1 & 1 & {\mathbf 0}^\top\\
 1 & 1 & 0 & (\widehat{{\mathbf a}}^{j-3})^\top\\
 1 & 0 & 1 & (\widehat{{\mathbf b}}^{j-3})^\top\\
 {\mathbf 0} & \widehat{{\mathbf a}}^{i-3} & \widehat{{\mathbf b}}^{i-3} & \widehat{C}_{j-3}^{i-3}
 \end {array}\right)
 $$
 $$
=\det\left(\begin {array} {cccc}
 1 & 0 & 0 & \bigl(\widehat{({\mathbf a}+{\mathbf b})}^{j-3}\bigr)^\top\\
 0 & 1 & 0 & (\widehat{{\mathbf b}}^{j-3})^\top\\
 0 & 0 & 1 & (\widehat{{\mathbf a}}^{j-3})^\top\\
\widehat{ ({\mathbf a}+ {\mathbf b})}^{i-3} & \widehat{{\mathbf
b}}^{i-3} & \widehat{{\mathbf a}}^{i-3} & \widehat{C}_{j-3}^{i-3}
 \end {array}\right)=b^{ij}_2,
 $$
 $$
b^{1j}_2=\det\left(\begin {array} {ccccc}
 0 & 1 & 0 & (\widehat{{\mathbf b}}^{j-3})^\top\\
 0 & 0 & 1 & (\widehat{{\mathbf a}}^{j-3})^\top\\
({\mathbf a}+ {\mathbf b}) & {\mathbf b}& {\mathbf a} &
\widehat{C}_{j-3}
 \end {array}\right)
 $$
 $$
=\det\left(\begin {array} {ccccc}
 1 & 1 & 0 & (\widehat{{\mathbf b}}^{j-3})^\top\\
 1 & 0 & 1 & (\widehat{{\mathbf a}}^{j-3})^\top\\
{\mathbf 0}& {\mathbf b}& {\mathbf a} & \widehat{C}_{j-3}
 \end {array}\right)=b^{1j}_1,
 $$
 $$
b^{2j}_2=\det\left(\begin {array} {cccc}
 1 & 0 & 0 & \bigl(\widehat{({\mathbf a}+{\mathbf b})}^{j-3}\bigr)^\top\\
 0 & 0 & 1 & (\widehat{{\mathbf a}}^{j-3})^\top\\
{\mathbf a}+ {\mathbf b} & {\mathbf b} & {\mathbf a} &
\widehat{C}_{j-3}
 \end {array}\right)
 $$
 $$
=\det\left(\begin {array} {cccc}
 1 & 0 & 0 & \bigl(\widehat{({\mathbf a}+{\mathbf b})}^{j-3}\bigr)^\top\\
 1 & 0 & 1 & (\widehat{{\mathbf a}}^{j-3})^\top\\
{\mathbf 0} & {\mathbf b} & {\mathbf a} & \widehat{C}_{j-3}
 \end {array}\right)
 =\det\left(\begin {array} {cccc}
 1 & 0 & 0 & (\widehat{{\mathbf b}}^{j-3})^\top\\
 1 & 0 & 1 & {\mathbf 0}^\top\\
{\mathbf 0} & {\mathbf b} & {\mathbf a} & \widehat{C}_{j-3}
 \end {array}\right)
 $$
 $$
 =\det\left(\begin {array} {cccc}
 1 & 1 & 0 & (\widehat{{\mathbf b}}^{j-3})^\top\\
 1 & 1 & 1 & {\mathbf 0}^\top\\
{\mathbf 0} & {\mathbf b} & {\mathbf a} & \widehat{C}_{j-3}
 \end {array}\right)=b^{2j}_1,
 $$
analogously $b^{3j}_1=b^{3j}_2$.

We have to verify that every vertex $x\notin\{v_1,\,v_2,\,v_3\}$ is
adjacent to either none of $v_1,\,v_2,\,v_3$ or precisely two of
them in $L_1$ and analogously for $L_2$. This statement is
equivalent to both $b^{1j}_1+b^{2j}_1+b^{3j}_1=0$ and
$b^{1j}_2+b^{2j}_2+b^{3j}_2=0$. By using the above equalities it is
enough to verify only the first equality.

We have
 $$
b^{2j}_1=\det\left(\begin {array} {cccc}
 1 & 1 & 1 & {\mathbf 0}^\top\\
 1 & 0 & 1 & (\widehat{{\mathbf b}}^{j-3})^\top\\
{\mathbf 0} & {\mathbf a} & {\mathbf b} & \widehat{C}_{j-3}
 \end {array}\right)=\det\left(\begin {array} {ccc}
  1 & 0 & (\widehat{{\mathbf b}}^{j-3})^\top\\
 {\mathbf a} & {\mathbf b} & \widehat{C}_{j-3}
 \end {array}\right),
 $$
 $$
b^{3j}_1=\det\left(\begin {array} {cccc}
 1 & 1 & 1 & {\mathbf 0}^\top\\
 1 & 1 & 0 & (\widehat{{\mathbf a}}^{j-3})^\top\\
{\mathbf 0} & {\mathbf a} & {\mathbf b} & \widehat{C}_{j-3}
 \end {array}\right)=\det\left(\begin {array} {ccc}
  0 & 1 & (\widehat{{\mathbf a}}^{j-3})^\top\\
 {\mathbf a} & {\mathbf b} & \widehat{C}_{j-3}
 \end {array}\right),
 $$
 $$
b^{1j}_1=\det\left(\begin {array} {ccccc}
 1 & 1 & 0 & (\widehat{{\mathbf a}}^{j-3})^\top\\
 1 & 0 & 1 & (\widehat{{\mathbf b}}^{j-3})^\top\\
{\mathbf 0}& {\mathbf a}& {\mathbf b} & \widehat{C}_{j-3}
 \end {array}\right)
 $$
 $$
=\det\left(\begin {array} {cccc}
 1 & 1 & \bigl(\widehat{({\mathbf a}+{\mathbf b})}^{j-3}\bigr)^\top\\
 {\mathbf a}& {\mathbf b} & \widehat{C}_{j-3}
 \end {array}\right)=b^{2j}_1+b^{3j}_1.
 $$

We have proven that our triples of vertices are related with each
other as in Remark~2.2, therefore, $L_1$ and $L_2$ are obtained from
each other by a sequence of $\Omega2$ and $\Omega3$ moves.
 \end {proof}

Let us define the map $\psi$ from the set of all homotopy classes of
looped graphs to the set of all graph-knots. Let $\mathfrak{L}$ be
the homotopy class of $L$. By using Lemma~\ref {lem:nondeg} we can
construct a symmetric matrix $\widetilde{A}(L)=(a_{ij})$ over
$\mathbb{Z}_2$ coinciding with the adjacency matrix of $L$ up to
diagonal elements and $\det \widetilde{A}(L)=1$. Let $G(L)$ be the
labeled simple graph having the matrix $\widetilde{A}(L)^{-1}+E$ as
its adjacency matrix (therefore, the first component of the label is
equal to the corresponding diagonal element of
$\widetilde{A}(L)^{-1}+E$), the second component of the label of the
vertex with the number $i$ is $w_i(1-2a_{ii})$, here $w_i=1$ if the
vertex of $L$ with the number $i$ does not have a loop, and $w_i=-1$
otherwise. Set $\psi(\mathfrak{L})=\mathfrak{F}$, here $G(L)$ is a
representative of $\mathfrak{F}$.

 \begin {theorem}\label {th:hom_grkn}
The map $\psi$ from the set of all homotopy classes of looped graphs
to the set of all graph-knots defined above is well-defined.
 \end {theorem}

 \begin {proof}
Let $L_1,\,L_2$ be two representatives of $\mathfrak{L}$. We have to
show that the graph-knots having representatives $G_1=G(L_1)$ and
$G_2=G(L_2)$ respectively are the same. We shall show that a
graph-link does not depend on the choice of diagonal elements, and
$G_1$ and $G_2$ are related to each other by Reidemeister
graph-moves.

$1)$ The independence (up to the second component of the label of a
vertex) on the choice of diagonal elements follows from the results
of~\cite {Ily}. Namely, graph-knots obtained from matrices having
different diagonal elements are related by $\Omega_g4$ and/or
$\Omega_g4'$ graph-moves (up to the second component of the label).
We have defined the sign of a vertex in such a way that looped
vertices correspond to vertices with the negative writhe number and
unlooped vertices correspond to vertices with the positive writhe
number (Lemma~\ref {lem:wrfr_si}). The writhe number doesn't change
under $\Omega_g4$, $\Omega_g4'$ graph-moves, and both the writhe
number and the framing allow one to determine the sign of a vertex.
Therefore, graph-knots obtained from matrices having different
diagonal elements are related by $\Omega_g4$ and/or $\Omega_g4'$
graph-moves.

Now we pass to moves on looped graphs.

$2)$ Let $L_1$ and $L_2$ be obtained from each other by $\Omega1$
(we remove the vertex with the number $1$). We have
 $$
A(L_1)=\left(\begin {array} {cc}
 a & {\mathbf 0}^\top\\
 {\mathbf 0} & A(L_2)
 \end {array}\right),\qquad a\in\{0,1\}.
 $$
Assume $\det\widetilde{A}(L_2)=1$, where $\widetilde{A}(L_2)$
coincides with $A(L_2)$ up to diagonal elements, and
 $$
\widetilde{A}(L_1)=\left(\begin {array} {cc}
 1 & {\mathbf 0}^\top\\
 {\mathbf 0} & \widetilde{A}(L_2)
 \end {array}\right).
 $$
Then
 $$
\widetilde{A}(L_1)^{-1}=\left(\begin {array} {cc}
 1 & {\mathbf 0}^\top\\
 {\mathbf 0} & \widetilde{A}(L_2)^{-1}
 \end {array}\right),
 $$
therefore, $G_1$ and $G_2$ are related by a sequence of $\Omega_g1$,
$\Omega_g4$, $\Omega_g4'$ graph-moves.

$3)$ Let $L_1$ and $L_2$ be obtained from each other by $\Omega2$
(we remove the vertex with the numbers $1$ and $2$). We have
 $$
A(L_1)=\left(\begin {array} {ccc}
 0 & b & {\mathbf a}^\top\\
 b & 1 & {\mathbf a}^\top\\
 {\mathbf a} & {\mathbf a}& A(L_2)
 \end {array}\right),\quad b\in\{0,1\}.
 $$
Assume $\det\widetilde{A}(L_2)=1$, where $\widetilde{A}(L_2)$
coincides with $A(L_2)$ up to diagonal elements, and
 $$
\widetilde{A}(L_1)=\left(\begin {array} {ccc}
 1+b & b & {\mathbf a}^\top\\
 b & 1+b & {\mathbf a}^\top\\
 {\mathbf a} & {\mathbf a}& \widetilde{A}(L_2)
 \end {array}\right),\quad\det\widetilde{A}(L_1)=1.
 $$
As
 $$
\det\left(\begin {array} {ccc}
 1+b & b & \widetilde{{\mathbf a}}^\top\\
 b & 1+b & \widetilde{{\mathbf a}}^\top\\
\widetilde{{\mathbf a}} & \widetilde{{\mathbf a}}& C
 \end {array}\right)=\det\left(\begin {array} {ccc}
 1+b & b & \widetilde{{\mathbf a}}^\top\\
 1 & 1 & {\mathbf 0}^\top\\
 \widetilde{{\mathbf a}} & \widetilde{{\mathbf a}} & C
 \end {array}\right)
 $$
 $$
=\det\left(\begin {array} {ccc}
 1+b & b & \widetilde{{\mathbf a}}^\top\\
 1 & 1 & {\mathbf 0}^\top\\
 {\mathbf 0} & {\mathbf 0}& C
 \end {array}\right)=\det C,
 $$
 $$
\det\left(\begin {array} {cc}
 1+b & {\mathbf a}^\top\\
 {\mathbf a}& \widetilde{A}(L_2)
 \end {array}\right)=\det\left(\begin {array} {cc}
 b & {\mathbf a}^\top\\
 {\mathbf a}& \widetilde{A}(L_2)
 \end {array}\right)+\det\widetilde{A}(L_2),
 $$
then
 $$
\widetilde{A}(L_1)^{-1}=\left(\begin {array} {ccc}
 f & f+1 & {\mathbf d}^\top\\
 f+1 & f & {\mathbf d}^\top\\
 {\mathbf d} & {\mathbf d}& \widetilde{A}(L_2)^{-1}
 \end {array}\right).
 $$
From the structure of the matrix $\widetilde{A}(L_1)^{-1}$ it
follows that the two vertices (which don't belong to $G_2$) have the
same framing and the necessary adjacencies, and from the structure
of the matrix $A(L_1)$ and the coincidence of vertices' framings it
follows that the two vertices have the different signs. This means
that $G_1$ and $G_2$ are obtained from each other by a sequence of
$\Omega_g2$, $\Omega_g4$, $\Omega_g4'$ graph-moves.

$4)$ Now assume that $L_1$ and $L_2$ are obtained from each other by
$\Omega 3$. Let us enumerate all the vertices of
$V(L_i)=\{v^i_1,\dots,v^i_n\}$ in such a way that corresponding
vertices of $L_1$ and $L_2$ under $\Omega 3$ move have the same
number, and without loss of generality we assume that $v^1_1$ and
$v^1_3$ are looped, $v^1_2$ is unlooped, $v^1_1$ is adjacent to
$v^1_2$; $v^1_3$ is adjacent to neither $v^1_1$ nor $v^1_2$ in
$L_1$. The case when $v^1_3$ is unlooped is obtained from the first
case by applying second Reidemeister and one of the third
Reidemeister moves.

We have
 $$
A(L_1)=\left(\begin {array} {cccc}
 1 & 1 & 0 & {\mathbf a}^\top\\
 1 & 0 & 0 & {\mathbf b}^\top\\
 0 & 0 & 1 & {\mathbf c}^\top\\
 {\mathbf a} & {\mathbf b} & {\mathbf c} & D
 \end {array}\right),\quad
A(L_2)=\left(\begin {array} {cccc}
 1 & 0 & 1 & {\mathbf a}^\top\\
 0 & 0 & 1 & {\mathbf b}^\top\\
 1 & 1 & 1 & {\mathbf c}^\top\\
 {\mathbf a} & {\mathbf b} & {\mathbf c} & D
 \end {array}\right),
 $$
 $$
{\mathbf a}+{\mathbf b}+{\mathbf c}= {\mathbf 0}.
 $$

Without loss of generality (if necessary we apply second
Reidemeister moves) we may assume that ${\mathbf c}\ne{\mathbf 0}$
and (Lemma~\ref {lem:nondeg})
 $$
\widetilde{A}(L_1)=\left(\begin {array} {cccc}
 0 & 1 & 0 & {\mathbf a}^\top\\
 1 & 1 & 0 & {\mathbf b}^\top\\
 0 & 0 & 0 & {\mathbf c}^\top\\
 {\mathbf a} & {\mathbf b} & {\mathbf c} & \widetilde{D}
 \end {array}\right)
 $$
and $\det\widetilde{A}(L_1)=1$.

As
 $$
\det\left(\begin {array} {cccc}
 0 & 0 & 1 & {\mathbf a}^\top\\
 0 & 0 & 1 & {\mathbf b}^\top\\
 1 & 1 & 1 & {\mathbf c}^\top\\
 {\mathbf a} & {\mathbf b} & {\mathbf c} & \widetilde{D}
 \end {array}\right)=\det\left(\begin {array} {cccc}
 0 & 0 & 1 & {\mathbf a}^\top\\
 0 & 0 & 0 & {\mathbf c}^\top\\
 1 & 0 & 1 &  {\mathbf c}^\top\\
 {\mathbf a} & {\mathbf c} &
{\mathbf c} & \widetilde{D} \end {array}\right)
 $$
 $$
=\det\left(\begin {array} {cccc}
 0 & 0 & 1 & {\mathbf a}^\top\\
 0 & 0 & 0 & {\mathbf c}^\top\\
 1 & 0 & 1 & {\mathbf b}^\top\\
 {\mathbf a} & {\mathbf c} &
{\mathbf b} & \widetilde{D}
 \end {array}\right)=\det\widetilde{A}(L_1)=1,\eqno(1)
 $$
we have
 $$
\widetilde{A}(L_2)=\left(\begin {array} {cccc}
 0 & 0 & 1 & {\mathbf a}^\top\\
 0 & 0 & 1 & {\mathbf b}^\top\\
 1 & 1 & 1 & {\mathbf c}^\top\\
 {\mathbf a} & {\mathbf b} & {\mathbf c} & \widetilde{D}
 \end {array}\right).
 $$
Let $\widetilde{A}(L_1)^{-1}=(\widetilde{a}^{ij}_1)$,
$\widetilde{A}(L_2)^{-1}=(\widetilde{a}^{ij}_2)$, and let
$G_i=G(L_i)$, $i=1,2$, be the two graph-knots having the adjacency
matrices $\widetilde{A}(L_i)^{-1}+E$. Let us show that $G_1$ and
$G_2$ are obtained from each other by a sequence of $\Omega_g 2$,
$\Omega_g 3$, $\Omega_g4$, $\Omega_g4'$ graph-moves.

Performing the same elementary manipulations as in $(1)$ we have
$\widetilde{a}^{ij}_1=\widetilde{a}^{ij}_2$ for $i,j\geqslant 3$.
Further, we get
 $$
1=\det\widetilde{A}(L_1)=\det\left(\begin {array} {cccc}
 0 & 1 & 0 & {\mathbf a}^\top\\
 1 & 1 & 0 & {\mathbf b}^\top\\
 0 & 0 & 0 & {\mathbf c}^\top\\
 {\mathbf a} & {\mathbf b} & {\mathbf c} & \widetilde{D}
 \end {array}\right)=\det\left(\begin {array} {cccc}
 0 & 1 & 0 & {\mathbf a}^\top\\
 1 & 1 & 0 & {\mathbf b}^\top\\
 1 & 0 & 0 & {\mathbf 0}^\top\\
 {\mathbf a} & {\mathbf b} & {\mathbf c} & \widetilde{D}
 \end {array}\right)
 $$
 $$
=\det\left(\begin {array} {ccc}
 1 & 0 & {\mathbf a}^\top\\
 1 & 0 & {\mathbf b}^\top\\
 {\mathbf b} & {\mathbf c} & \widetilde{D}
 \end {array}\right)=\det\left(\begin {array} {ccc}
 1 & 0 & {\mathbf b}^\top\\
 0 & 0 & {\mathbf c}^\top\\
 {\mathbf b} & {\mathbf c} & \widetilde{D}
 \end {array}\right)=\widetilde{a}^{11}_1,
 $$
 $$
\widetilde{a}^{12}_1=\det\left(\begin {array} {ccc}
 1 & 0 & {\mathbf b}^\top\\
 0 & 0 & {\mathbf c}^\top\\
 {\mathbf a} & {\mathbf c} & \widetilde{D}
 \end {array}\right)=\det\left(\begin {array} {ccc}
 1 & 0 & {\mathbf b}^\top\\
 0 & 0 & {\mathbf c}^\top\\
 {\mathbf b} & {\mathbf c} & \widetilde{D}
 \end {array}\right)=1,
 $$
 $$
\widetilde{a}^{13}_1=\det\left(\begin {array} {ccc}
 1 & 1 & {\mathbf b}^\top\\
 0 & 0 & {\mathbf c}^\top\\
 {\mathbf a} & {\mathbf b} & \widetilde{D}
 \end {array}\right)=\det\left(\begin {array} {ccc}
 0 & 1 & {\mathbf b}^\top\\
 0 & 0 & {\mathbf c}^\top\\
 {\mathbf c} & {\mathbf b} & \widetilde{D}
 \end {array}\right)=1,
 $$
 $$
\widetilde{a}^{1j}_1=\det\left(\begin {array} {cccc}
 1 & 1 & 0 & (\widehat{{\mathbf b}}^{j-3})^\top\\
 0 & 0 & 0 & (\widehat{{\mathbf c}}^{j-3})^\top\\
{\mathbf a} & {\mathbf b} & {\mathbf c} &
\widehat{\widetilde{D}}_{j-3}
 \end {array}\right)=0,\quad j\geqslant 3\quad({\mathbf a}+{\mathbf b}+{\mathbf c}= {\mathbf
 0}),
 $$
 $$
\widetilde{a}^{2j}_1=\det\left(\begin {array} {cccc}
 0 & 1 & 0 & (\widehat{{\mathbf a}}^{j-3})^\top\\
 0 & 0 & 0 & (\widehat{{\mathbf c}}^{j-3})^\top\\
{\mathbf a} & {\mathbf b} & {\mathbf c} &
\widehat{\widetilde{D}}_{j-3}
 \end {array}\right)=\det\left(\begin {array} {cccc}
 0 & 1 & 0 & (\widehat{{\mathbf a}}^{j-3})^\top\\
 0 & 0 & 0 & (\widehat{{\mathbf c}}^{j-3})^\top\\
{\mathbf a} & {\mathbf 0} & {\mathbf c} &
\widehat{\widetilde{D}}_{j-3}
 \end {array}\right),
 $$
 $$
\widetilde{a}^{3j}_1=\det\left(\begin {array} {cccc}
 0 & 1 & 0 & (\widehat{{\mathbf a}}^{j-3})^\top\\
 1 & 1 & 0 & (\widehat{{\mathbf b}}^{j-3})^\top\\
{\mathbf a} & {\mathbf b} & {\mathbf c} &
\widehat{\widetilde{D}}_{j-3}
 \end {array}\right)=\det\left(\begin {array} {cccc}
 0 & 1 & 0 & (\widehat{{\mathbf a}}^{j-3})^\top\\
 1 & 0 & 0 & (\widehat{{\mathbf b}}^{j-3})^\top\\
{\mathbf a} & {\mathbf 0} & {\mathbf c} &
\widehat{\widetilde{D}}_{j-3}
 \end {array}\right).
 $$
If either $\widetilde{a}^{22}_1=0$ or $\widetilde{a}^{33}_1=0$ we
can apply the same second Reidemeister graph-moves to $G_1$ and
$G_2$ and after that applying the $\Omega'_g4$ graph-move we get
that the corresponding vertices have the framing $0$. Analogously,
if $\widetilde{a}^{32}_1=1$ we can apply the same second
Reidemeister graph-moves to $G_1$ and $G_2$ and after that applying
the $\Omega_g4$ graph-move we get that $v_2^1$ and $v^1_3$ are not
adjacent to each other.  Therefore, without loss of generality we
may assume that $\widetilde{a}^{22}_1=\widetilde{a}^{33}_1=1$,
$\widetilde{a}^{32}_1=0$. By using the last equalities we get
 $$
\widetilde{a}^{22}_1=\det\left(\begin {array} {ccc}
 0 & 0 & {\mathbf a}^\top\\
 0 & 0 & {\mathbf c}^\top\\
 {\mathbf a} & {\mathbf c} & \widetilde{D}
 \end {array}\right)=\det\left(\begin {array} {ccc}
 0 & 0 & {\mathbf b}^\top\\
 0 & 0 & {\mathbf c}^\top\\
 {\mathbf b} & {\mathbf c} & \widetilde{D}
 \end {array}\right)
 $$
 $$
 =1+\det\left(\begin {array} {cc}
  0 & {\mathbf c}^\top\\
 {\mathbf c} & \widetilde{D}
 \end {array}\right)=1,
 $$
 $$
\widetilde{a}^{33}_1=\det\left(\begin {array} {ccc}
 0 & 1 & {\mathbf a}^\top\\
 1 & 1 & {\mathbf b}^\top\\
 {\mathbf a} & {\mathbf b} & \widetilde{D}
 \end {array}\right)=\det\left(\begin {array} {ccc}
 1 & 0 & {\mathbf c}^\top\\
 0 & 1 & {\mathbf b}^\top\\
 {\mathbf c} & {\mathbf b} & \widetilde{D}
 \end {array}\right)
 $$
 $$
 =1+\det\left(\begin {array} {cc}
  1 & {\mathbf b}^\top\\
 {\mathbf b} & \widetilde{D}
 \end {array}\right)=1,
 $$
 $$
\widetilde{a}^{23}_1=\det\left(\begin {array} {ccc}
 0 & 1 & {\mathbf a}^\top\\
 0 & 0 & {\mathbf c}^\top\\
 {\mathbf a} & {\mathbf b} & \widetilde{D}
 \end {array}\right)=\det\left(\begin {array} {ccc}
 1 & 1 & {\mathbf b}^\top\\
 0 & 0 & {\mathbf c}^\top\\
 {\mathbf c} & {\mathbf b} & \widetilde{D}
 \end {array}\right)
 $$
 $$
 =1+\det\left(\begin {array} {cc}
  0 & {\mathbf c}^\top\\
 {\mathbf b} & \widetilde{D}
 \end {array}\right)=0.
 $$

Let us find the remaining elements of $\widetilde{A}(L_2)^{-1}$. We
have
 $$
\widetilde{a}^{11}_2=\det\left(\begin {array} {ccc}
 0 & 1 & {\mathbf b}^\top\\
 1 & 1 & {\mathbf c}^\top\\
 {\mathbf b} & {\mathbf c} & \widetilde{D}
 \end {array}\right)=
 \det\left(\begin {array} {ccc}
 0 & 1 & {\mathbf b}^\top\\
 1 & 0 & {\mathbf a}^\top\\
 {\mathbf b} & {\mathbf c} & \widetilde{D}
 \end {array}\right)=
 \det\left(\begin {array} {ccc}
 1 & 1 & {\mathbf b}^\top\\
 1 & 0 & {\mathbf a}^\top\\
 {\mathbf b} & {\mathbf c} & \widetilde{D}
 \end {array}\right)
 $$
 $$
 +\det\left(\begin {array} {cc}
 0 & {\mathbf a}^\top\\
 {\mathbf c} & \widetilde{D}
 \end {array}\right)=
 \det\left(\begin {array} {ccc}
 1 & 1 & {\mathbf b}^\top\\
 0 & 1 & {\mathbf c}^\top\\
 {\mathbf b} & {\mathbf c} & \widetilde{D}
 \end {array}\right)+\det\left(\begin {array} {cc}
 0 & {\mathbf b}^\top\\
 {\mathbf c} & \widetilde{D}
 \end {array}\right)+\det\left(\begin {array} {cc}
 0 & {\mathbf c}^\top\\
 {\mathbf c} & \widetilde{D}
 \end {array}\right)
 $$
 $$
=
 \det\left(\begin {array} {ccc}
 1 & 1 & {\mathbf b}^\top\\
 0 & 0 & {\mathbf c}^\top\\
 {\mathbf b} & {\mathbf c} & \widetilde{D}
 \end {array}\right)+\det\left(\begin {array} {cc}
 1 & {\mathbf b}^\top\\
 {\mathbf b} & \widetilde{D}
 \end {array}\right)+1=1,
 $$
 $$
\widetilde{a}^{22}_2=\det\left(\begin {array} {ccc}
 0 & 1 & {\mathbf a}^\top\\
 1 & 1 & {\mathbf c}^\top\\
 {\mathbf a} & {\mathbf c} & \widetilde{D}
 \end {array}\right)=
 \det\left(\begin {array} {ccc}
 1 & 1 & {\mathbf b}^\top\\
 1 & 1 & {\mathbf c}^\top\\
 {\mathbf b} & {\mathbf c} & \widetilde{D}
 \end {array}\right)=1,
 $$
 $$
\widetilde{a}^{33}_2=\det\left(\begin {array} {ccc}
 0 & 0 & {\mathbf a}^\top\\
 0 & 0 & {\mathbf b}^\top\\
 {\mathbf a} & {\mathbf b} & \widetilde{D}
 \end {array}\right)=
 \det\left(\begin {array} {ccc}
 0 & 0 & {\mathbf b}^\top\\
 0 & 0 & {\mathbf c}^\top\\
 {\mathbf b} & {\mathbf c} & \widetilde{D}
 \end {array}\right)=1,
 $$
 $$
\widetilde{a}^{12}_2=\det\left(\begin {array} {ccc}
 0 & 1 & {\mathbf b}^\top\\
 1 & 1 & {\mathbf c}^\top\\
 {\mathbf a} & {\mathbf c} & \widetilde{D}
 \end {array}\right)=
\det\left(\begin {array} {ccc}
 1 & 1 & {\mathbf b}^\top\\
 0 & 1 & {\mathbf c}^\top\\
 {\mathbf b} & {\mathbf c} & \widetilde{D}
 \end {array}\right)
 $$
 $$
=\det\left(\begin {array} {ccc}
 1 & 1 & {\mathbf b}^\top\\
 0 & 0 & {\mathbf c}^\top\\
 {\mathbf b} & {\mathbf c} & \widetilde{D}
 \end {array}\right)+\det\left(\begin {array} {cc}
 1 & {\mathbf b}^\top\\
 {\mathbf b} & \widetilde{D}
 \end {array}\right)=0,
 $$
 $$
\widetilde{a}^{13}_2=\det\left(\begin {array} {ccc}
 0 & 0 & {\mathbf b}^\top\\
 1 & 1 & {\mathbf c}^\top\\
 {\mathbf a} & {\mathbf b} & \widetilde{D}
 \end {array}\right)=
\det\left(\begin {array} {ccc}
 0 & 0 & {\mathbf b}^\top\\
 0 & 1 & {\mathbf c}^\top\\
 {\mathbf c} & {\mathbf b} & \widetilde{D}
 \end {array}\right)
 $$
 $$
= \det\left(\begin {array} {ccc}
 0 & 0 & {\mathbf b}^\top\\
 0 & 0 & {\mathbf c}^\top\\
 {\mathbf c} & {\mathbf b} & \widetilde{D}
 \end {array}\right)+\det\left(\begin {array} {cc}
 0 & {\mathbf b}^\top\\
 {\mathbf c} & \widetilde{D}
 \end {array}\right)=0,
 $$
 $$
\widetilde{a}^{23}_2=\det\left(\begin {array} {ccc}
 0 & 0 & {\mathbf a}^\top\\
 1 & 1 & {\mathbf c}^\top\\
 {\mathbf a} & {\mathbf b} & \widetilde{D}
 \end {array}\right)=
\det\left(\begin {array} {ccc}
 0 & 1 & {\mathbf b}^\top\\
 0 & 1 & {\mathbf c}^\top\\
 {\mathbf c} & {\mathbf b} & \widetilde{D}
 \end {array}\right)=0,
 $$
 $$
\widetilde{a}^{1j}_2=\det\left(\begin {array} {cccc}
 0 & 0 & 1 & (\widehat{{\mathbf b}}^{j-3})^\top\\
 1 & 1 & 1 & (\widehat{{\mathbf c}}^{j-3})^\top\\
{\mathbf a} & {\mathbf b} & {\mathbf c} &
\widehat{\widetilde{D}}_{j-3}
 \end {array}\right)=\det\left(\begin {array} {cccc}
 0 & 0 & 1 & (\widehat{{\mathbf b}}^{j-3})^\top\\
 1 & 1 & 0 & (\widehat{{\mathbf a}}^{j-3})^\top\\
{\mathbf a} & {\mathbf b} & {\mathbf c} &
\widehat{\widetilde{D}}_{j-3}
 \end {array}\right)
 $$
 $$
=\det\left(\begin {array} {cccc}
 1 & 0 & 1 & (\widehat{{\mathbf b}}^{j-3})^\top\\
 0 & 1 & 0 & (\widehat{{\mathbf a}}^{j-3})^\top\\
{\mathbf 0} & {\mathbf b} & {\mathbf c} &
\widehat{\widetilde{D}}_{j-3}
 \end {array}\right)=\det\left(\begin {array} {ccc}
1 & 0 & (\widehat{{\mathbf a}}^{j-3})^\top\\
{\mathbf b} & {\mathbf c} & \widehat{\widetilde{D}}_{j-3}
 \end {array}\right)=\widetilde{a}^{2j}_1+\widetilde{a}^{3j}_1,
 $$
 $$
\widetilde{a}^{2j}_2=\det\left(\begin {array} {cccc}
 0 & 0 & 1 & (\widehat{{\mathbf a}}^{j-3})^\top\\
 1 & 1 & 1 & (\widehat{{\mathbf c}}^{j-3})^\top\\
{\mathbf a} & {\mathbf b} & {\mathbf c} &
\widehat{\widetilde{D}}_{j-3}
 \end {array}\right)=\det\left(\begin {array} {cccc}
 0 & 1 & 1 & (\widehat{{\mathbf a}}^{j-3})^\top\\
 1 & 0 & 0 & (\widehat{{\mathbf b}}^{j-3})^\top\\
{\mathbf a} & {\mathbf 0} & {\mathbf c} &
\widehat{\widetilde{D}}_{j-3}
 \end {array}\right)=\widetilde{a}^{3j}_1,
 $$
 $$
\widetilde{a}^{3j}_2=\det\left(\begin {array} {cccc}
 0 & 0 & 1 & (\widehat{{\mathbf a}}^{j-3})^\top\\
 0 & 0 & 1 & (\widehat{{\mathbf b}}^{j-3})^\top\\
{\mathbf a} & {\mathbf b} & {\mathbf c} &
\widehat{\widetilde{D}}_{j-3}
 \end {array}\right)=\det\left(\begin {array} {cccc}
 0 & 0 & 1 & (\widehat{{\mathbf a}}^{j-3})^\top\\
 0 & 0 & 0 & (\widehat{{\mathbf c}}^{j-3})^\top\\
{\mathbf a} & {\mathbf b} & {\mathbf 0} &
\widehat{\widetilde{D}}_{j-3}
 \end {array}\right)=\widetilde{a}^{2j}_1.
 $$

We see that
 $$
\widetilde{A}(L_1)^{-1}+E=\left(\begin {array} {cccc}
 0 & 1 & 1 & {\mathbf 0}^\top\\
 1 & 0 & 0 & {\mathbf f}^\top\\
 1 & 0 & 0 & {\mathbf g}^\top\\
 {\mathbf 0} & {\mathbf f} & {\mathbf g} & H
 \end {array}\right),
 $$
 $$
\widetilde{A}(L_2)^{-1}+E=\left(\begin {array} {cccc}
 0 & 0 & 0 & {\mathbf f}^\top+{\mathbf g}^\top\\
 0 & 0 & 0 & {\mathbf g}^\top\\
 0 & 0 & 0 & {\mathbf f}^\top\\
 {\mathbf f}+{\mathbf g} & {\mathbf g} & {\mathbf f} & H
 \end {array}\right).
 $$
It is easy to see that the corresponding vertices have the structure
as in Definition~2.9. Therefore, $G_2$ is obtained from $G_1$ by a
sequence of $\Omega_g 2$, $\Omega_g 3$, $\Omega_g4$, $\Omega_g4'$
graph-moves.
 \end {proof}

From Theorems~\ref {th:grkn_ho},~\ref {th:hom_grkn} and from the
definitions of the map $\chi$ and $\psi$ it follows that these maps
are mutually inverse. Therefore, we have proven Theorem~\ref
{th:gk_hl}.

We conclude the paper with the example of non-realisable graph-knot.

 \begin {definition}
A graph-link (a homotopy class of looped graphs) is called {\em
non-realisable} if it has no realisable representative.
 \end {definition}

 \begin {crl}
A graph-link $\mathfrak{F}$ is non-realisable if and only if
$\chi(\mathfrak{F})$ is non-realisable.
 \end {crl}

\begin{figure} \centering\includegraphics[width=100pt]{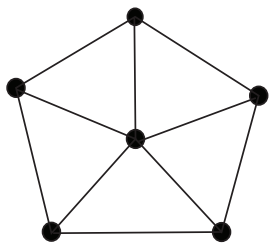}
 \caption{The First Bouchet Graph Gives a Non-Realisable
Graph-Knot} \label{Bouch}
 \end{figure}

Let $G$ be the labeled graph depicted in Fig.~\ref {Bouch} with each
vertex having the framing $1$ and signs are chosen arbitrary.

 \begin {crl}
The graph-knot $\mathfrak{F}$ generated by $G$ is non-realisable.
 \end {crl}

 \begin {proof}
Let $\mathfrak{L}=\chi(\mathfrak{F})$. It is not difficult to see
that the looped graph isomorphic to $G$ (we have no loop) is a
representative of $\mathfrak{L}$. Therefore, $\mathfrak{L}$ is
non-realisable, see~\cite {Ma4,Ma5}, and, in turn, $\mathfrak{F}$ is
non-realisable graph-knot.
 \end {proof}

 \end {document}